\documentclass[10pt]{birkjour}
 \newtheorem{thm}{Theorem}[section]
 \newtheorem{cor}[thm]{Corollary}
 \newtheorem{lem}[thm]{Lemma}
 \newtheorem{prop}[thm]{Proposition}
 \theoremstyle{definition}
 \newtheorem{defn}[thm]{Definition}
 \theoremstyle{remark}
 \newtheorem{rem}[thm]{Remark}
 \newtheorem*{ex}{Example}
 \numberwithin{equation}{section}
 \usepackage{mathges}
 \band{38}{2018}
 \author{Helmut Kahl}
 \usepackage{amsmath}
 \usepackage{amssymb,textcomp}
 \usepackage{graphicx}
 \usepackage{epstopdf}
\title{Measuring quadric sectors at centre
\thanks{The author thanks to H.-J. Kroll for careful reading and useful hints.}
}

\begin{document}
\maketitle

\centerline{Helmut Kahl}
\vglue4pt
{\bf Abstract.}
Sectors at centre of quadrics with point symmetry are investigated over arbitrary fields of characteristic different from two. As an application nice formulas are demonstrated for the area and the volume of such planar and spatial sectors, respectively, in euclidean space. It seems that up to now there has been at most little research in this field up to very special cases.

{\bf Mathematical Subject Cassification(2010).} Primary 51N10; Secondary 28A75.


\section{Introduction}
Quadrics represent the most simple non-linear algebraic varieties. In dimension two and three they were investigated already in the Greek-Hellenistic antiquity (see \cite{Scriba}, sect.2.2.2, p.42 and sect.2.5.10, p.92). Sectors of such quadrics centred at a point of symmetry (not to be confused with Kepler{\textquotesingle}s sectors at a focus of an ellipse in his {\textquotedblleft}Astronomia Nova{\textquotedblright}; see \cite{Scriba}, sect.5.2, p.266-268) seem to have been of little interest up to now. Only special cases have been considered intensively; e.g.: a circle sector is essential for the notion of an angle, a sector at the symmetry point of a unit hyperbola is essential for interpreting the inverse hyperbolic functions $\operatorname{arcosh} , \operatorname{arsinh}$, and spherical triangles have been important to astronomy for about two thousand years (see \cite{Scriba}, sect.2.5.6, p.81). Nevertheless also those general sectors at centre deserve attention for the following reason: Because of its simple geometric nature they allow
\begin{itemize}
\item a simple arithmetic description and thus might be useful for {\textquotedblleft}second order approximations{\textquotedblright} of arbitrary curves and surfaces,
\item a simple area and volume formula (see Figure 1 and Theorems \ref{thm_area}, \ref{thm_volume}) and thus might be useful for {\textquotedblleft}second order approximations{\textquotedblright} of arbitrary areas and volumes.
\end{itemize}
In section 2 an elementary theory of quadric sectors at centre in vector spaces of finite dimension over commutative fields of characteristic different from two will be developed. In section 3 its measure theory over the field $\mathbb{R}$ of real numbers will be represented with general formulas in $\mathbb{R}^{2}$ and $\mathbb{R}^{3}$. The key to the formulas is the so called \textit{sector coefficient}, a certain affine invariant. The additivity (see Corollary \ref{cor_angle}) of the generalised angle of Definition \ref{def_angle} emphasises the usefulness of the sector coefficient for measuring quadric sectors. It is desirable to find an integral free formula in case of a hyperboloid of one and of two sheets. It seems to the author that up to now there has been at most little research in the field of quadric sectors at centre, although their quadrics are represented elegantly by the sector coefficients due to Proposition \ref{prop_sector}.

\section{Affine quadric sectors at centre}
In this section $\mathbb{K}$ denotes a commutative field of characteristic different from two. For a finite-dimensional vector space $V$ over $\mathbb{K}$ we denote by $\mathbb{A}(V , \mathbb{K})$ the affine space with point set $V$ and with the cosets $v + U$ of one-dimensional subspaces $U$ of $V$ as lines ($v \in V$).\footnote{Every affine plane fulfilling the axiom of Desargues and every at least 3-dimensional affine space can be represented this way (see \cite{KSW}, Satz(10.1)\&(11.20)).} The zero vector is denoted by $o$. For every translation $\tau$ of $\mathbb{A}(V , \mathbb{K})$ there exists a vector $t \in V$ s.t. $\tau(x) = x + t$ for every $x \in V$ (see \cite{KSW}, Satz(12.2)). Every automorphism $\alpha$ of $\mathbb{A}(V , \mathbb{K})$ is a composition $\alpha = \tau \circ \iota$ of a semilinear bijection $\iota$ of $V$ with a translation $\tau$ (see \cite{KSW}, Satz(12.6)). If $\iota$ is linear $\alpha$ is called a (\textit{linear}) \textit{affinity}. For a subset $X \subseteq V$ denote $\langle X \rangle$ and $\overline{X}$ the linear hull and incidence hull of $X$, respectively. A point $c \in V$ is called a \textit{centre} of $X$ if $2 c - x \in X$ for all $x \in X$.

\begin{lem}\label{lem_translation}
The set $C$ of all centres of a non-empty set $X \subseteq V$ fulfills:


a) For $c \in C , x \in X$ and a subspace $A \supseteq X$ holds $c = \frac{1}{2} (2 c - x + x) \in A$.

b) For all $c \in C$ holds $\overline{X} = c + \langle X - c \rangle$.

c) If $X$ is a subspace then $C = X$.

d) If $C \ne \emptyset$ is a subspace and $\tau$ a translation with $\tau(C) = C$ then $\tau(X) = X$.
\end{lem}
\begin{proof}
a) For $c \in C , x \in X$ and a subspace $A \supseteq X$ holds $c = \frac{1}{2} (2 c - x + x) \in A$.

b) This follows from a).

c) i) By a) holds $C \subseteq X$. ii) For $c , x \in X$ holds $2 c - x \in X$ since $X$ is a subspace. This shows $c \in C$, hence $X \subseteq C$.

d) For $c \in C$ holds $c' := \tau(c) \in C$, hence $c'' := (c + c') / 2 \in C$. So for every $x \in X$ holds $\tau(x) = x + c' - c = 2 c'' - \left(2 c - x \right) \in X$.
\end{proof}

A non-empty set $Q \subseteq V$ is called a \textit{quadric} of $V$ when there is a quadratic form $q:V \to \mathbb{K}$, a linear form $l:V \to \mathbb{K}$ and a scalar $\gamma \in \mathbb{K}$ with
\begin{equation}\label{quadric}
Q = \lbrace x \in V \mid q(x)+l(x)+\gamma=0 \rbrace .
\end{equation}
For every affinity $\Phi : V \to V'$ the image $\Phi(Q)$ is a quadric of $V'$.

\begin{ex}
For an affinity $\Phi : \mathbb{R}^3 \to \mathbb{R}^3$ the quadric
\begin{equation*}
Q := \Phi\left(\lbrace (x,y,z) \in \mathbb{R}^3 \mid x^2 + y^2 = 1 \rbrace\right)
\end{equation*}
is an \textit{elliptic cylinder}. With '$-$' instead of '$+$' it is a \textit{hyperbolic cylinder}. In both cases the set $C$ of centres is a line.
\end{ex}

\begin{rem}\label{rem_quadric}
Let $q , l , \gamma$ define $Q$ by equation \eqref{quadric}.

a) For a point $c \in V$ there is a linear from $l_c$ and a scalar $\gamma_c$ s.t.
\begin{equation*}
Q = \lbrace x \in V \mid q(x-c) + l_c(x-c) + \gamma_c = 0 \rbrace .
\end{equation*}
Since the \textit{polar form}
\begin{equation*}
\varphi(x,y):= \left(q(x+y) - q(x) - q(y)\right)/2
\end{equation*}
of $q$ is bilinear and fulfills $q(x-c) = q(x) + q(c) - 2 \varphi(x,c)$ this assertion holds for $l_c(x) := l(x) + 2 \varphi(x,c)$ and $\gamma_c := q(c)+l(c)+\gamma$.

b) If $c$ is a centre of the quadric $Q_c := \lbrace x \in V \mid q(x-c) + l(x-c) + \gamma = 0 \rbrace$ we have $l(x-c)=0$ and $q(x-c) = -\gamma$ for all $x \in Q_c$. This follows from the two equations that arise from $x , 2 c - x \in Q_c$.

c) In case $Q = \lbrace o \rbrace$ the only centre of $Q$ is $o$. Otherwise let $\lbrace b_1 , ... , b_n \rbrace$ be a basis of $\langle Q \rangle$. Then the number of centres of $Q$ is one if and only if the $n \times n$-matrix $(\varphi(b_i , b_j))$ is invertible. Proof: An $x \in V$ is a centre of $Q$ if and only if for all $p \in Q$ holds
\begin{equation*}
4 q(x) + q(p) - 4 \varphi(x,p) + 2 l(x) - l(p) + \gamma = 0 .
\end{equation*}
Because of $q(p)+l(p)+\gamma=0$ that equation is equivalent with
\begin{equation*}
2 q(x) + l(x) = 2 \varphi(x,p) + l(p) .
\end{equation*}
Since the latter term is linear in $p$ it follows: $x$ is an element of the set $C$ of centres if and only if $2 \varphi(x,p) + l(p) = 0$ for all $p \in \langle Q \rangle$. But for the linear forms $\tilde{\varphi}_x(y) := \varphi(x,y)$ the equation $\tilde{\varphi}_x = -l /2$ over $\langle Q \rangle$ has a unique solution $x \in \langle Q \rangle$ if and only if $(\varphi(b_i , b_j))$ is invertible. Because of $C \subseteq \overline{Q} \subseteq \langle Q \rangle$ this shows the assertion.
\end{rem}

\begin{lem}\label{lem_centre}
The set $C$ of centres of a quadric $Q$ fulfills:

a) $C$ is a subspace.

b) If $C \cap Q$ is not empty then $Q = \underset{x \in Q}{\bigcup} \overline{C \cup \lbrace x \rbrace} \supseteq C$.

c) For $q , l , \gamma$ defining $Q$ by equation \eqref{quadric} and a centre $c \in C \setminus Q$ it holds
\begin{equation*}
Q = \lbrace x \in V \mid q(x-c) + \gamma_c = 0 \rbrace \textnormal{ with } \gamma_c := q(c)+l(c)+\gamma \ne 0 .
\end{equation*}
\end{lem}
\begin{proof}
a) Let $C$ be non-empty. Via translation we may assume $o \in C$ without loss of generality. Then $x \in Q = \lbrace x \in V \mid q(x)+l(x)+\gamma=0 \rbrace$ implies $l(x) = 0$ by Remark \ref{rem_quadric}b). So with the polar form $\varphi$ of $q$ it holds for $x \in Q$ and $c \in C$
\begin{equation*}
-4 \varphi(c , x) = 2 \varphi(2 c , -x) = q(2 c - x) - q(2 c) - q(x) = - q(2 c) = -4 q(c) ,
\end{equation*}
therefore $\varphi(c , x) = q(c)$. Since $x \mapsto \varphi(c , x)$ is linear and because $o \in \overline{Q}$ due to Lemma \ref{lem_translation}a) it follows $\varphi(c , x) = \varphi(c , 0) = 0$ for all $x \in \langle Q \rangle = \overline{Q}$. For the subspaces $H_x := \lbrace y \in \langle Q \rangle \mid \varphi(y , x) = 0 \rbrace$ that means \begin{equation*}
C \subseteq C' := \underset{x \in \langle Q \rangle}{\bigcap} H_x .
\end{equation*}
For $c \in C' \subseteq \langle Q \rangle$ holds $q(c) = \varphi(c , c) = 0$. So for $x \in Q$ we have
\begin{equation*}
q(2 c - x) - q(x) = q(2 c - x) - q(2 c) - q(x) = 2 \varphi(2 c , -x) = -4 \varphi(c , x) = 0 .
\end{equation*}
Therefore $c$ is a centre which proofs $C = C'$. So $C$ is a subspace.

b) By Remark \ref{rem_quadric}a)\&b) for $c \in C \cap Q$ there is a quadratic form $q$ and a linear form $l$ s.t. $Q = \lbrace x \in V \mid q(x-c) + l(x-c) = 0 \rbrace$ and $q(x-c) = l(x-c) = 0$ for all $x \in Q$. Hence for $c' \in C$ we have $0 = q(2 c ' - c - c) = 4 q(c' - c)$ and $0 = l(2 c ' - c - c) = 2 l(c' - c)$, thus $c' \in Q$. Therefore it holds $C \subseteq Q$. So for arbitrary $c \in C$ we have the $Q$ defining equation above. For $x \in Q$ and $\lambda \in \mathbb{K}$ it follows $q(c + \lambda (x - c) - c) + l(c + \lambda (x - c) - c) = \lambda^2 q(x - c) + \lambda(x - c)  = 0$. This shows that the subspace $\overline{c , x} := c + \langle x - c \rangle$ is contained in $Q$. Because of Lemma \ref{lem_translation}d) the subspace $x + \langle C - c \rangle$ through $x \in Q$ which is parallel $C$ is also contained in $Q$. Therefore it holds
\begin{equation*}
\overline{C \cup \lbrace x \rbrace} = \left(x + \langle C - c \rangle\right) \cup \underset{c \in C}{\bigcup} \overline{c , x} \subseteq Q ,
\end{equation*}
hence the equation.

c) Via translation we may assume $c = o$ without loss of generality. By Remark \ref{rem_quadric}b) $Q \subseteq H := \lbrace x \in V \mid l(x) = 0 \rbrace$ and $\gamma \ne 0$ because of $c \notin Q$. Hence there is some $a \in Q$ with $\alpha := q(a) \ne 0$, since $Q$ is not empty. Assume that $l$ does not vanish. Then there is some $b \in V \setminus H$ with $\varphi(a,b) = 0$, since for arbitrary $b \in V \setminus H$ the vector
\begin{equation*}
b - \frac{\varphi(a,b)}{\alpha} a
\end{equation*}
is still not in $H$. By definition of $H$ it holds $l(b) \ne 0$. Assume $q(b) \ne 0$. Then
\begin{equation*}
c := a - \frac{l(b)}{q(b)} b \in V \setminus H
\end{equation*}
would fulfill $q(c)+l(c) = q(a) + l(a)$, in contradiction to $Q \subseteq H$. So we have
\begin{equation*}
q\left(-\frac{\gamma}{l(b)} b\right) + l\left(-\frac{\gamma}{l(b)} b\right) + \gamma = 0 - \gamma + \gamma = 0 ,
\end{equation*}
also in contradiction to $Q \subseteq H$. So $l$ must vanish.
\end{proof}

\begin{ex}
For an affinity $\Phi : \mathbb{K}^3 \to \mathbb{K}^3$ the \textit{cone}
\begin{equation*}
Q := \Phi\left(\lbrace (x,y,z) \in \mathbb{K}^3 \mid x^2 + y^2 - z^2 = 0 \rbrace\right)
\end{equation*}
possesses the centre $c := \Phi(0 , 0 , 0) \in Q$ only. The union of all lines $\overline{c , x}$ with $x \in Q \setminus \lbrace c \rbrace$ gives $Q$. The common intersection point of those lines is $c$.
\end{ex}

From now on we restrict to quadrics $Q$ with an \textit{external centre} $p_0 \in C$, i.e. $p_0 \in C \setminus Q$. So by Lemma \ref{lem_centre}b) it holds $C \cap Q = \emptyset$. By Lemma \ref{lem_centre}c) there is a quadratic form $q$ s.t.
\begin{equation*}
Q = \lbrace x \in V \mid q(x-p_0) = 1 \rbrace .
\end{equation*}
Since $Q$ is not empty, there is a point $p_1 \in Q$. Assume another point $p_2 \in Q$ with $p_2 \ne p_1$ and $2 p_0 - p_1 \ne p_2$. Then $p_1 - p_0$ and $p_2 - p_0$ are linearly independent because for $\lambda \in \mathbb{K}$ the equation $1 = q(\lambda (p_1 - p_0)) = \lambda^2$ implies $\lambda = 1$ or $\lambda = -1$.

\begin{defn}\label{def_sector}
For a quadric $Q$ with external centre $p_0$ and points $p_1 , ... , p_n \in Q$ s.t. $p_1 - p_0 , ... , p_n - p_0$ are linearly independent we call
\begin{equation*}
\sigma := \left(Q , p_0 ; p_1 , ... , p_n\right)
\end{equation*}
an $n$-\textit{dimensional sector at centre} $p_0$ \textit{of} $Q$ \textit{with vertices} $p_{1},p_{2},...,p_{n}$. The function $\Phi_\sigma :\mathbb{K}^{n}\to V$ defined by
\begin{equation*}
\Phi_\sigma(x_1 , ... , x_n) := p_0 + \sum\limits_{i=1}^{n} x_i (p_i - p_0)
\end{equation*}
is called its \textit{frame affinity}. A two- or three-dimensional sector is called \textit{planar} or \textit{spatial}, respectively.
\end{defn}

\begin{rem}
Let $\sigma := \left(Q , p_0 ; p_1 , ... , p_n\right)$ denote a sector of a quadric $Q$ and let $Q_\sigma := Q \cap \Phi_\sigma \left(\mathbb{K}^{n}\right)$.

a) If $\Phi_\sigma(\mathbb{K}^{n})$ is a vector space its subset $Q_\sigma$ is a quadric.

b) In case $\dim V = n$ the set $Q_\sigma$ coincides with $Q$.
\end{rem}

\begin{prop}\label{prop_sector}
For a sector $\sigma = \left(Q , p_0 ; p_1 , ... , p_n\right)$ let $\varphi$ denote the polar form of a quadratic form $q$ s.t. $Q = \lbrace x \in V \mid q(x-p_0) = 1 \rbrace$ and let
\begin{equation*}
\sigma _{i j} := \varphi(p_i - p_0 , p_j - p_0)
\end{equation*}
for all $i , j \in \{1,...,n\}$. Then $S := \left(\sigma _{i j}\right)$ is the unique symmetric $n \times n$-matrix  over $\mathbb{K}$ s.t. $\Phi_\sigma^{-1} (Q) = \lbrace x \in \mathbb{K}^n \mid x S x^t = 1 \rbrace$. For all $i \in \{1,...,n\}$ holds $\sigma _{i i} = 1$.
\end{prop}
\begin{proof}
The vectors $b_i := p_i - p_0$ fulfill
\begin{equation*}
\varphi(x_1 b_1 + ... + x_n b_n , y_1 b_1 + ... + y_n b_n) = x S y^t
\end{equation*}
for all $x = (x_1 , ... , x_n) , y = (y_1 , ... , y_n) \in \mathbb{K}^n$. This implies $q(\Phi_\sigma(x)-p_0) = x S x^t$, hence the claimed identity for $S$. Since
\begin{equation*}
e_1:=(1,0,...,0), e_2:=(0,1,0,...,0), ... , e_n:=(0,...,0,1)
\end{equation*}
are elements of $P := \Phi_\sigma^{-1} (Q)$ we have $\sigma_{i i} = 1$. It remains to prove uniqueness.\footnote{see also \cite{AuHa}, Satz 2} Let $(\tilde{\sigma}_{i j})$ be another symmetric matrix that defines $P$. Then $\tilde{\sigma}_{i i} = 1$ follows as above. So it suffices to show $\tilde{\sigma}_{i j} = \sigma_{i j}$ for $i < j$. First assume $n = 2$ and let $\beta := 2 \sigma_{1 2}$. That means $x_1^2 + \beta x_1 x_2 + x_2^2 = (x_1,x_2) S (x_1,x_2)^t$ for all $x_1 , x_2 \in \mathbb{K}$. Let $\tilde{\beta} := 2 \tilde{\sigma}_{1 2}$. Since $(\beta,-1) \in P$ we get $\beta^2 - \beta \tilde{\beta} + 1 = 1$, hence $\beta^2 = \beta \tilde{\beta}$. Analogously, it follows $\tilde{\beta}^2 = \beta \tilde{\beta}$. These identities imply $\tilde{\beta} = \beta$, ie. $\tilde{\sigma}_{1 2} = \sigma_{1 2}$. In case $n > 2$, for given $i , j$ with $i < j$, this argumentation applies to the point $(x_1 , ... , x_n) \in P$ defined by $x_i := 2 \sigma_{i j} , x_j := -1$ and $x_k := 0$ for all $k \in \{1,...,n\} \setminus \lbrace i , j \rbrace$.
\end{proof}

\begin{defn}\label{def_sectorCoeff}
Due to Proposition \ref{prop_sector} for a planar sector $\sigma = (Q , p_0 ; p_1 , p_2)$ the field element $\chi(\sigma):=\beta$ is well-defined by
\begin{equation*}
\Phi^{-1}_\sigma(Q) = \lbrace (x,y) \in \mathbb{K}^2 \mid x^{2} + \beta x y + y^{2} = 1 \rbrace .
\end{equation*}
It is called the \textit{sector coefficient of} $\sigma$.
\end{defn}

\begin{ex}
For every injective linear map $\Phi: \mathbb{K}^{n}\to W$ into a vector space $W$ the image $Q := \Phi(U)$ of the quadric $U := \lbrace x \in \mathbb{K}^{n} \mid x x^t = 1\rbrace$ is a quadric of the linear subspace $V := \Phi(\mathbb{K}^{n})$ of $W$. With canonical unit vectors $e_i \in U$ we obtain an $n$-dimensional sector $\left(Q, o ; \Phi(e_1) , ... , \Phi(e_n)\right)$ that fulfills $\chi\left(Q, o ; \Phi(e_i) ,\Phi(e_j)\right) = 0$ for all $i \ne j$.
\end{ex}

\begin{rem}\label{rem_sectorCoeff}
a) By Proposition \ref{prop_sector} the sector coefficient is invariant under linear affinities.

b) For $x , y \in \mathbb{K} \setminus \lbrace 0 \rbrace$ with $\Phi_\sigma(x,y) \in Q$ it holds
\begin{equation*}
\chi(\sigma) = \frac{1}{x y}-\frac{x}{y}-\frac{y}{x} .
\end{equation*}

c) Let $Q \subset \mathbb{K}^2$ be a quadric with external centre $p_0 := (0 , 0)$ and with pairwise linearly independent $p_1 := (\alpha , \beta) , p_2 := (\gamma , \delta) , p_3 := (\varepsilon , \zeta) \in Q$. Then for
\begin{equation*}
\lambda := \left|\begin{matrix} \alpha & \beta \\ \gamma & \delta \end{matrix}\right| , \mu := \left|\begin{matrix} \varepsilon & \zeta \\ \gamma & \delta \end{matrix}\right| , \nu := \left|\begin{matrix} \alpha & \beta \\ \varepsilon & \zeta \end{matrix}\right|
\end{equation*}
it holds $p_3 = \frac{\mu}{\lambda} p_1 + \frac{\nu}{\lambda} p_2$, hence by Remark b):
\begin{equation*}
\chi(Q,p_0;p_1,p_2) = \frac{\lambda^2}{\mu \nu} - \frac{\mu}{\nu} - \frac{\nu}{\mu}
\end{equation*}
\end{rem}

d) In Remark c) the rational function
\begin{equation*}
r(x , y , z) := \left( z^2 - x^2 - y^2 \right) / (2 x y)
\end{equation*}
is used for describing a sector coefficient. It has another algebraic property: For $\lambda, \mu, \nu \in \mathbb{K}\setminus\lbrace 0 \rbrace$ let $\delta := r(\lambda, \mu, \nu) , \varepsilon := r(\mu, \nu, \lambda) , \zeta := r(\nu, \lambda, \mu)$. Then it holds
\begin{equation*}
\delta^2 + \varepsilon^2 + \zeta^2 - 2 \delta \varepsilon \zeta = 1 .
\end{equation*}

\begin{thm}\label{thm_sectorCoeff}
a) Let $\sigma := (Q , p_0 ; p_1 , ... , p_n)$ and $\sigma' := (Q' , p_0 ; p_1 , ... , p_n)$ be sectors of dimension $n>1$. When for all $i , j \in \lbrace 1 , ... , n \rbrace$ with $i < j$ holds $Q \cap \overline{p_0 , p_i , p_j} = Q' \cap \overline{p_0 , p_i , p_j}$ then $Q \cap \Phi_\sigma(\mathbb{K}^n) = Q' \cap \Phi_{\sigma'}(\mathbb{K}^n)$.\footnote{I.e.: Those $n(n-1)/2$ planar intersections with $Q$ determine $Q \cap \Phi_\sigma(\mathbb{K}^n)$.}

b) Let $\lbrace b_1 , ... , b_n \rbrace$ be a basis of $V$. For all $i , j \in \lbrace 1 , ... , n \rbrace$ with $i < j$ let $x_i , y_j \in \mathbb{K} \setminus \lbrace 0 \rbrace$. Then there is only one quadric with external centre $o$ that contains all $b_i$ and all $x_i b_i + y_j b_j$.

c) Let $V$ have finite dimension. If $Q$ is a quadric with external centre $c \in V$ and $\overline{Q}=V$ then there is a unique quadratic form $q$ s.t. $Q = \lbrace x \in V \mid q(x-c)=1 \rbrace$. Changing the centre $c$ changes $q$ only by a non-zero factor.
\end{thm}
\begin{proof}
a) By definition of the frame affinity it holds $\Phi_\sigma = \Phi_{\sigma'}$. Via translation we may assume $p_0 = o$. Then $Q_{i j} := Q \cap \overline{p_0 , p_i , p_j} = Q \cap \langle p_i , p_j \rangle$ is a quadric of $\langle p_i , p_j \rangle$. Hence for the symmetric matrix $S = (\sigma_{i j})$ that corresponds by Proposition \ref{prop_sector} with $\sigma$ we have $\chi(Q_{i j} , o ; p_i , p_j) = 2 \sigma_{i j}$ for $i < j$. So $S$ corresponds also with $\sigma'$, and it follows $\Phi_\sigma^{-1}(Q') = \Phi_\sigma^{-1}(Q)$.

b) For all $i < j$ set $\sigma_{j i} := \sigma_{i j} := \frac{1}{2}\left(\frac{1}{x_i y_j}-\frac{x_i}{y_j}-\frac{y_j}{x_i}\right)$, and $\sigma_{i i} := 1$ for all $i$. Then $S := (\sigma_{i j})$ is symmetric. For $P := \lbrace x \in \mathbb{K}^n \mid x S x^t = 1 \rbrace$ and $\Phi(x_1 , ... , x_n) := x_1 b_1 + ...+ x_n b_n$ the quadric $Q := \Phi(P)$ contains the $b_i$. But it contains also the $x_i b_i + y_j b_j$ since for the canonical unit vectors $e_i \in \mathbb{K}^n$ holds
\begin{equation*}
(x_i e_i + y_j e_j) S (x_i e_i + y_j e_j)^t = x_i^2 + y_j^2 + 2 \sigma_{i j} x_i y_j = 1 .
\end{equation*}
For another quadric $Q'$ with that property we have $\chi(Q' , o ; b_i , b_j) = 2 \sigma_{i j}$ by Remark \ref{rem_sectorCoeff}b). So by Proposition \ref{prop_sector} it follows $Q' = \Phi(S) = Q$ since $\Phi$ is the frame affinity of $(Q , o ; b_1 , ... , b_n)$ and of $(Q' , o ; b_1 , ... , b_n)$.

c) Let $n := \dim V$ be the dimension of $V$. Because of $\overline{Q}=V$ and Lemma \ref{lem_translation}b) there is an $n$-dimensional sector $\sigma = (Q , p_0 ; p_1 , ... , p_n)$ with $p_0 := c$. With its symmetric matrix $S$ of Proposition \ref{prop_sector} we define $q$ by $q\left(\Phi_\sigma(x) - p_0 \right) =  x S x^t$ for $x \in \mathbb{K}^n$. Then the equation $Q = \lbrace y \in V \mid q(y-p_0)=1 \rbrace$ holds. That shows the existence.\footnote{For uniqueness see also \cite{AuHa}, Satz 2 again.} For the polar form $\varphi$ of $q$ and $x_1 , ... , x_n \in \mathbb{K}$ we have
\begin{equation*}
q \left(\Phi_\sigma(x_1 , ... , x_n) - p_0 \right) = \sum\limits_{i=1}^{n} x_i^2 + 2 \sum\limits_{i<j} \varphi(p_i - p_0 , p_j - p_0) x_i x_j ,
\end{equation*}
whereby the second sum is taken over all pairs $(i , j) \in \{1 , ... , n\} \times \{1 , ... , n\}$ with $i < j$. Because of uniqueness of $S$ it follows $2 \varphi(p_i - p_0 , p_j -p_0) = \chi(Q , p_0; p_i , p_j)$. Since $\varphi$ is bilinear and $p_1 - p_0 , ... , p_n - p_0$ confirm a basis $\varphi$ is determined by those sector coefficients, hence $q$ is also. Now let $c'$ be another centre of $Q$. By Lemma \ref{lem_centre}c) it holds $Q = \lbrace x \in V \mid q(x - c') = \gamma \rbrace$ for some $\gamma \ne 0$. Therefrom the second assertion follows.
\end{proof}

\begin{prop}\label{prop_fullSector}
Let $V$ be a vector space of finite dimension $n$ over a field with more than five elements and with $1+1 \ne 0$. Then a quadric $Q \subset V$ with external centre possesses a sector of dimension $n$, hence it holds $\overline{Q} = V$.
\end{prop}
\begin{proof}
Via translation we may assume that $o$ is a centre of $Q$. Then by Lemma \ref{lem_centre}c) it holds $Q = \lbrace x \in V | q(x) = 1 \rbrace$ for some quadratic form $q$. We proof the assertion by induction on $n$. For $n=1$ it is clear. For $n > 1$ let $H$ be a linear hyperspace of $V$ with some $p \in Q \cap H$. Choose $b \in V \setminus H$. Then for $a := b - \varphi(b,p) p$ holds $\varphi(a,p)=0$, hence $q(\lambda p + \mu a) = \lambda^2 + q(a) \mu^2$ for $\lambda , \mu \in \mathbb{K}$. But for $\alpha := q(a)$ there is a $\beta \ne 0$ s.t $\alpha \beta^2 \ne 1$ and $\alpha \beta^2 \ne -1$, because $\mathbb{K}$ has more than five elements. Hence there are $\lambda , \mu \ne 0$ s.t. $\lambda^2 + \alpha \mu^2 = 1$, namely
\begin{equation*}
\lambda := \frac{\alpha \beta^2 - 1}{\alpha \beta^2 + 1} , \mu := \frac{2 \beta}{\alpha \beta^2 + 1} .
\end{equation*}
Therefore we have $p_n := \lambda p + \mu a \in Q \setminus H$. By induction hypothesis there are linearly independent $p_1 , ... , p_{n-1} \in Q \cap H$. So $(Q,o;p_1 , ... , p_n)$ is a sector.
\end{proof}

\begin{cor}\label{cor_sector}
Let $V$ be a finite-dimensional vector space over a field with more than five elements and with $1+1 \ne 0$. Then for every $c \in V$ the map
\begin{equation*}
q \mapsto \lbrace x \in V \mid q(x-c) = 1 \rbrace .
\end{equation*}
from the set of all quadratic forms $q$ of $V$ that allow a solution $x \in V$ of the equation $q(x)=1$ to the set of all quadrics of $V$ with external centre $c$ is bijective.
\end{cor}
\begin{proof}
This follows from Theorem \ref{thm_sectorCoeff}c) and Proposition \ref{prop_fullSector}.
\end{proof}

\section{Measure of a sector at centre in euclidean space}
For a quadric $Q$ of the euclidean space $\mathbb{R}^n$ with external centre $p_0$ we consider sectors $\sigma := (Q , p_0 ; p_1 , ... , p_n)$ of full dimension $n$. We may assume $p_0 = o = (0,...,0)$. Let us call the symmetric $n \times n$-matrix $S$ of Proposition \ref{prop_sector} the \textit{sector matrix} of $\sigma$. For an $x \in \mathbb{R}^n$ we write $x \geqslant 0$ when every coordinate of $x$ is non-negative. Then for $P^{+} := \left\{x \in \mathbb{R}^{n} \mid x \geqslant 0 , \; x S x^t \leqslant 1 \right\}$ and the frame affinity $\Phi(x_1, ... ,x_n) := x_1 p_{1} + ... + x_n p_{n}$ we define the \textit{sector region} to be measured as the set
\begin{equation*}
Q^{+} := \Phi\left(P^{+}\right) .
\end{equation*}
In case $n=2$ or $n=3$ the set $P^{+}$ is the region in the main quadrant or octant, respectively, bounded by $P := \Phi^{-1}(Q)$. In general $Q^{+}$ possesses the vertices $o,p_{1},...,p_{n}$. It is bounded by the $n(n-1)/2$ planes $\left\langle p_{i},p_{j}\right\rangle \subset \mathbb{R}^{n}$ ($1\leqslant i < j \leqslant n$) and by $Q$. But it is not necessarily \textit{bounded} in the following sense: There is a constant bounding the vector norm of $x$ for all $x\in Q^{+}$ . We will measure only bounded sector regions. Therefore we exclude the quadrics that contain a centre since they are unions of certain subspaces (cf. Lemma \ref{lem_centre}b)).

\textit{Examples}.
a) When $p_1 , p_2$ lie on two different branches of a hyperbola or on two different parallel lines, the corresponding planar sector region is not bounded.

b) A spherical triangle with vertices $p_{1},p_{2},p_{3}$ together with the sphere{\textquotesingle}s centre $o$ determines a bounded sector region in $\mathbb{R}^{3}$.

c) Let $Q \subset \mathbb{R}^{3}$ be a circular cylinder and $p_{1},p_{2},p_{3}\in Q$ three pairwise different points. If the line $Z$ of $Q$'s centres intersects the triangle between $p_{1},p_{2},p_{3}$ in a point different from $o$ then the sector region of $\left(Q,o;p_{1},p_{2},p_{3}\right)$ is not bounded since it contains {\textquotedblleft}half{\textquotedblright} of the line $Z$.

For the volume measure function $\mu$ of $\mathbb{R}^{n}$ and a linear function $\Psi : \mathbb{R}^{n} \to \mathbb{R}^{n}$ we use the transformation formula (see \cite{Bronstein}, no. 8.147)
\begin{equation}\label{transfFormula}
\mu(\Psi(M)) = |\det \Psi| \mu(M)
\end{equation}
for all measurable $M \subset \mathbb{R}^{n}$. So we have $\mu \left(Q^{+}\right) = |\det \Phi| \mu\left(P^{+}\right)$. If we denote by $\Theta_n$ the volume of the $n$-simplex with vertices $o,p_{1},...,p_{n}$, we get
\begin{equation}\label{volume}
\mu \left(Q^{+}\right) = \Theta_n \: n! \: \mu\left(P^{+}\right) .
\end{equation}
For $n=2$ and $n=3$ the factor $\Theta_n$ is well known; e.g. in dependence of the {\textquotedblleft}geodesic data{\textquotedblright} length $|a|$ of a vector $a$ and measure of an angle:
\begin{itemize}
\item The area of a triangle between edges $a,b$ with their angle of measure $\omega \in (0,\pi )$ is $\Theta_2 = \frac{1}{2} |a| |b| \sin(\omega)$.
\item The volume of a tetrahedron between edges $a,b,c$ with angles of measure $\varphi ,\psi ,\omega \in (0,\pi )$ between these edges is
\begin{equation*}
\Theta_3 = \frac{1}{6} |a| |b| |c| \sqrt{ \left|\det \left(\begin{matrix}1&\cos \varphi &\cos \psi \\\cos \varphi &1&\cos \omega \\\cos \psi &\cos \omega &1\end{matrix}\right)\right| } \;.
\end{equation*}
\end{itemize}

The following remark shows how to compute a sector coefficient in dependence of lengths and angles.

\begin{rem}\label{rem_geodesy}
The sector coefficient of a planar sector $(Q,o;a,b)$ is determined by a $c \in Q$ that is linearly independent from $a$ and from $b$. For the measure $\varphi$ of the \textit{oriented angle} (see \cite{Audin}, chap.III.1) from $a$ to $c$ and the measure $\psi$ of the oriented angle from $c$ to $b$ we get by some trigonometry
\begin{equation*}
x = \frac{|c| \sin \psi}{|a| \sin (\varphi +\psi)} \; , \; y = \frac{|c| \sin \varphi}{|b| \sin (\varphi +\psi)}
\end{equation*}
when $c = x a + y b$. Hence the questionable sector coefficient reads
\begin{equation}\label{geodata}
\chi(Q,o;a,b) = \frac{|a| |b| \sin^{2}(\varphi +\psi)}{|c|^{2}\sin \varphi \sin \psi}-\frac{|b| \sin \psi}{|a| \sin \varphi}-\frac{|a| \sin \varphi}{|b| \sin \psi}\;.
\end{equation}
In the special case of a circle it equals $2 \cos \omega$, where $\omega :=\varphi +\psi $ denotes the measure of the angle between sides $a$ and $b$ (oriented or non-oriented - it doesn{\textquotesingle}t matter). By choosing $\psi =\varphi $ this follows from the double-angle formula $2\sin ^{2}\varphi =1-\cos (2\varphi )$.
\end{rem}

The factor $\mu\left(P^{+}\right)$ of formula \eqref{volume} is determined by the sector coefficients $2 \sigma_{i j}$ ($1\leqslant i<j\leqslant n$) of the given sector. The main task is to express it as an analytic function of the $n(n-1)/2$ variables $\sigma_{i j}$. From now we restrict to the most practical cases $n=2$ and $n=3$. Hereby we use the inverse trigonometric functions $\arccos, \arcsin$ and the inverse hyperbolic functions $\operatorname{arcosh}, \operatorname{arsinh}$.

\begin{center}
\includegraphics[width=11.6cm,height=8.2cm]{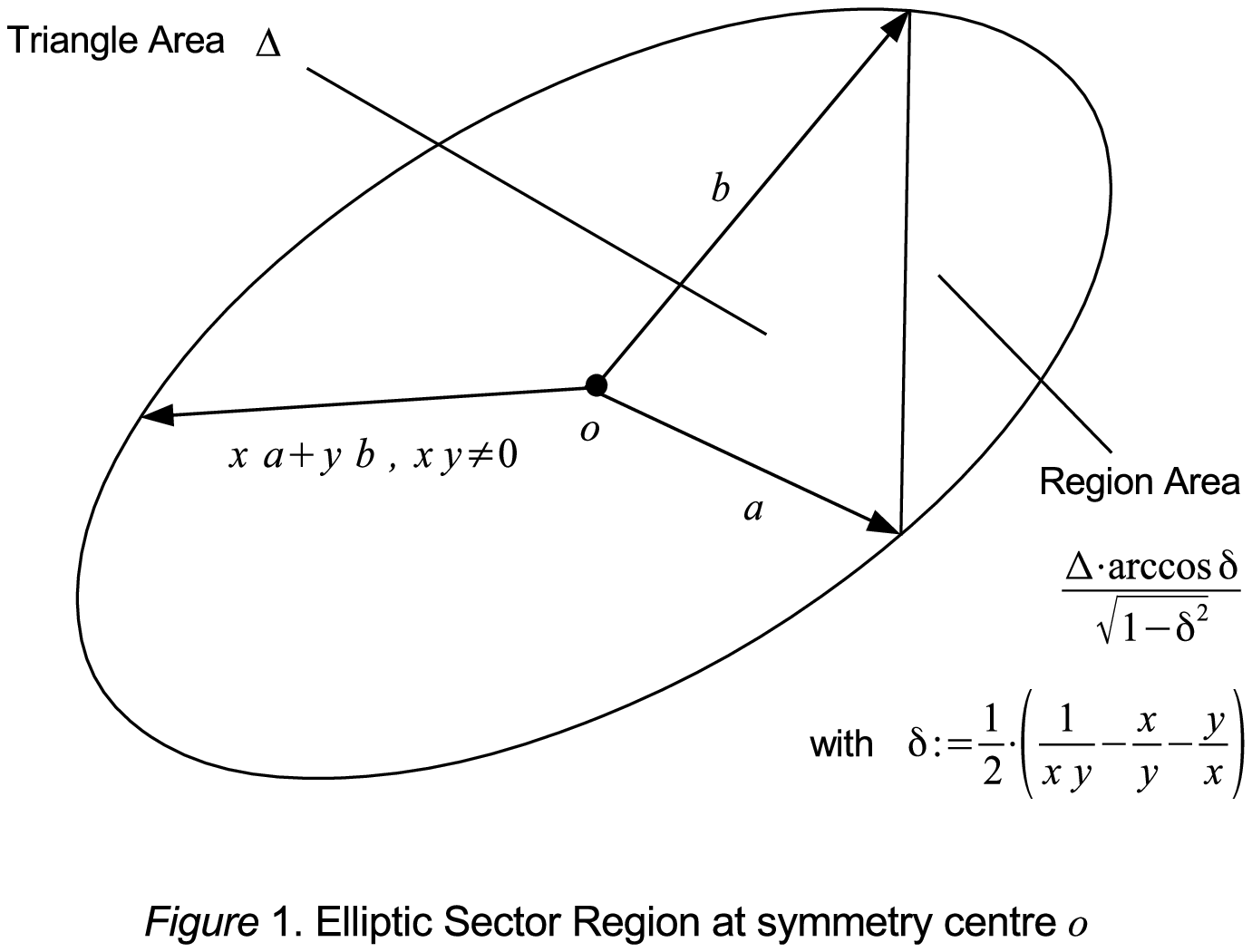}
\end{center}

\begin{thm}\label{thm_area}
The region of a planar sector $(Q,o;a,b)$ is bounded if and only if $\delta := \chi(Q,o;a,b) / 2 > -1$. Then its area $\mu(Q^{+})$ amounts to
\begin{equation*}
\mu(Q^{+}) = \left\lbrace \begin{matrix} \Delta {\frac{\arccos\delta}{\sqrt{1-\delta ^{2}}}} & \textnormal{in the elliptic case} \hfill\null & |\delta| < 1 \\ \Delta \hfill\null & \textnormal{in the straight line case} \hfill\null & \delta =1 \\ \Delta {\frac{\operatorname{arcosh}\delta }{\sqrt{\delta ^{2}-1}}} & \textnormal{in the hyperbolic case} \hfill\null & \delta > 1 \end{matrix}\right.
\end{equation*}
with $\Delta := \Theta_2$ denoting the area of the triangle with vertices $o,a,b$.
\end{thm}
\begin{proof}
According to Proposition \ref{prop_sector} the region in question is the image of
\begin{equation*}
P^{+} = \lbrace (x,y) \in \mathbb{R}^2 \mid x^{2} + 2 \delta x y + y^{2} \leqslant 1; x , y \geqslant 0 \rbrace
\end{equation*}
under $\Phi(x,y) := x a + y b$. Hence it is bounded if and only if $P^{+}$ is bounded. Since $P^{+}$ is symmetric with respect to the bisecting line $y=x$, this is the case if and only if $P^{+}$ contains a non-zero point on this line, i.e.  $\delta >-1$. In case $|\delta |<1$ we substitute
\begin{equation*}
\left(\begin{matrix}x\\y\end{matrix}\right) \text{ by } \left(\begin{matrix} 1 & \delta \\ 0 & \sqrt{1-\delta ^{2}}\end{matrix}\right) \left(\begin{matrix}x\\y\end{matrix}\right)\;,
\end{equation*}
thus getting the unit circle $U := \lbrace (x,y) \in \mathbb{R}^2 \mid x^{2}+y^{2} = 1 \rbrace$ from $P = \lbrace (x,y) \in \mathbb{R}^2 \mid x^{2}+2\delta xy+y^{2} = 1 \rbrace$. The area of the circle sector $U^{+}$ between $(1,0)$ and $\left(\delta ,\sqrt{1-\delta ^{2}}\right)$ is $(\arccos \delta )/2$. Now the claimed formula in the elliptic case follows by formula \eqref{transfFormula}, since the functional determinant of the transformation from $U^{+}$ to $Q^{+}$ is $\Delta \sqrt{1-\delta^{2}}$. The case $\delta =1$ is trivial. In case $\delta >1$ we substitute
\begin{equation*}
\left(\begin{matrix}x\\y\end{matrix}\right) \text{ by } \left(\begin{matrix} 1 & \delta \\ 0 & \sqrt{\delta ^{2}-1}\end{matrix}\right) \left(\begin{matrix}x\\y\end{matrix}\right)\;,
\end{equation*}
thus getting the unit hyperbola $U := \lbrace (x,y) \in \mathbb{R}^2 \mid x^{2}-y^{2} = 1 \rbrace$. The area of the hyperbola sector $U^{+}$ between $(1,0)$ and $\left(\delta ,\sqrt{\delta ^{2}-1}\right)$ is $(\operatorname{arcosh}\delta )/2$. Now the claimed formula in the hyperbolic case follows in the same way like in the elliptic case.
\end{proof}

\begin{rem}\label{rem_determinants}
The proof of Theorem \ref{thm_area} has shown that any sector coefficient determines the type of the planar quadric. So by Theorem \ref{thm_sectorCoeff}c) and Remark \ref{rem_sectorCoeff}c) we obtain the following property of determinants: Let $a_1 , a_2 , b_1 , b_2 , c_1 , c_2$ real numbers s.t. $\lambda := a_1 b_2 - a_2 b_1 , \mu := b_1 c_2 - b_2 c_1 , \nu := c_1 a_2 - c_2 a_1$ do not vanish; and let $\delta := r(\lambda, \mu, \nu) , \varepsilon := r(\mu, \nu, \lambda) , \zeta := r(\nu, \lambda, \mu)$ with
\begin{equation*}
r(x , y , z) := \left( z^2 - x^2 - y^2 \right) / (2 x y) .
\end{equation*}
Then all $\kappa \in \lbrace \delta , \varepsilon , \zeta \rbrace$ fulfill simultaneously $|\kappa| < 1$ or $|\kappa| = 1$ or $|\kappa| > 1$.
\end{rem}

\begin{thm}\label{thm_volume}
For a spatial sector $\sigma := (Q,o;p_1,p_2,p_3)$ let $\Theta := \Theta_3$ be the volume of its tetrahedron with vertices $o,p_1,p_2,p_2$. For the $3\times 3$-sector matrix $S=(\sigma _{i j})$ of $\sigma$ set $\eta :=\sigma_{1 2} , \theta :=\sigma_{2 3} , \kappa :=\sigma_{1 3}$. We may assume $\eta \leqslant \theta \leqslant \kappa $ by suitable permutation of the vertices $p_1,p_2,p_3$. Then the volume of the sector region $Q^{+}$ of $\sigma$ is $\mu \left(Q^{+}\right) = \Theta f(\eta ,\theta ,\kappa)$ for the analytic function $f$ whose term $f(\eta ,\theta ,\kappa)$ is described in the following:
\begin{equation*}
\frac{2}{\sqrt{|S|}} \left(a(\eta,\theta,\kappa) + a(\kappa,\eta,\theta) + a(\theta,\kappa,\eta) - \pi \right)
\end{equation*}
with
\begin{equation*}
a(\eta,\theta,\kappa) := \arccos \frac{\eta -\theta \kappa }{\sqrt{\left(1-\theta ^{2}\right)\left(1-\kappa ^{2}\right)}}
\end{equation*}
in case $|S|=1+2\eta \theta \kappa -\eta ^{2}-\theta ^{2}-\kappa ^{2}>0 \wedge |\eta |,|\theta |,|\kappa |<1$ of an ellipsoid;
\begin{equation*}
\frac{2}{1+\eta } \text{ or } 2 (1-\eta ) \frac{\kappa +\theta -\eta -1}{(\theta -\eta \kappa )(\kappa -\eta \theta )}
\end{equation*}
in case $|S|=0 \wedge -1<\eta \wedge (\kappa = 1 > \theta $ or $\theta >\eta \kappa)$, respectively, of an elliptic cylinder;
\begin{equation*}
f = 1 \text{ in case } (\eta ,\theta ,\kappa) = (1,1,1) \text{ of a plane;}
\end{equation*}
\begin{equation*}
\frac{2}{1+\kappa } \text{ or } 2 (\kappa -1) \frac{1+\kappa -\theta -\eta }{(\eta \kappa -\theta )(\theta \kappa -\eta )}
\end{equation*}
in case $|S|=0 \wedge (\eta = 1 <\theta $ or $\eta >1)$, respectively, of a hyperbolic cylinder;

$\frac{2}{\sqrt{\kappa^{2} - 1}} \int_{0}^{1}\left(\operatorname{arsinh}\frac{\left(\kappa^{2} - 1\right) \left(\sqrt{1 + \eta^{2}y^{2} - y^{2}}-\eta y\right)+(\kappa \theta - \eta)y}{\sqrt{\kappa^{2}-1+|S|y^{2}}} - \operatorname{arsinh}\frac{(\kappa \theta - \eta)y}{\sqrt{\kappa^{2} - 1 + |S|y^{2}}}\right)\mathit{dy}$

in case $|S|>0 \wedge 1<\eta$ of a hyperboloid of two sheets and in case $|S|<0 \wedge 1\leqslant \eta$ of a hyperboloid of one sheet;

$\frac{2}{\sqrt{1 - \eta^{2}}} \int_{0}^{1}\left(\arcsin \frac{\left(1 - \eta^{2}\right) \left(\sqrt{1 + \kappa^{2}y^{2} - y^{2}} - \kappa y \right)+(\kappa - \eta \theta)y}{\sqrt{1 - \eta^{2} - |S|y^{2}}} - \arcsin \frac{(\kappa - \eta \theta)y}{\sqrt{1 - \eta^{2} - |S|y^{2}}}\right)\mathit{dy}$

in case $|S|<0 \wedge 0 \leqslant \eta < 1$ of a hyperboloid of one sheet.

The sector region $Q^{+}$ is also bounded for $|S|<0 , -1 < \eta < 0 , \theta > \eta \kappa$. This case of a hyperboloid of one sheet is reduced to the case $|S|<0$ with all sector coefficients being positive by splitting $\sigma$ into the four sectors
\begin{equation*}
(Q,o;p_1,p_{1 2},p_{1 3}), (Q,o;p_{1 2},p_2,p_{2 3}), (Q,o;p_{1 3},p_{2 3},p_3), (Q,o;p_{1 2},p_{2 3},p_{1 3})
\end{equation*}
with
\begin{equation*}
p_{j k} := \frac{p_j + p_k}{\sqrt{2 + \chi(Q,o;p_j,p_k)}} \text{ for } j < k .
\end{equation*}
All other cases of $\eta ,\theta ,\kappa $ which are not described for the above six types of quadrics do not yield bounded sector regions.
\end{thm}

\begin{proof}
The sector matrix $S$ defines the quadric
\begin{equation*}
P := \lbrace (x,y,z) \in \mathbb{R}^3 \mid x^{2} + 2 \eta x y + y^{2} + 2 \theta y z + z^{2} + 2 \kappa z x = 1 \rbrace .
\end{equation*}
Due to Theorem \ref{thm_area} the inequalities $|\eta |,|\theta |,|\kappa |<1$ must hold in case of an ellipsoid. We transform $P$ by the isomorphism $\Psi(x,y,z) := T \cdot (x,y,z)^t$ with
\begin{equation*}
T:=\left(\begin{matrix} 1 & \eta & \kappa \\ 0 & \sqrt{1-\eta ^{2}} & (\theta -\eta \kappa)/\sqrt{1-\eta ^{2}} \\ 0 & 0 & \sqrt{|S|} / \sqrt{1-\eta ^{2}}\end{matrix}\right)\;,
\end{equation*}
thus getting the unit sphere $\Psi(P) = \lbrace (x,y,z) \in \mathbb{R}^3 \mid x^{2} + y^{2} + z^{2} = 1 \rbrace$. Under isomorphisms the sign of the determinant of the quadric defining matrix does not change. Hence we must have indeed $|S|>0$ in the ellipsoid case. Because of $|T|=\sqrt{|S|}$ the volume of the spherical triangle $\Psi\left(P^{+}\right)$ is $\mu(\Psi\left(P^{+}\right)) = \mu (P^{+}) \sqrt{|S|}$ due to formula \eqref{transfFormula}. Its edge vectors are the columns of $T$. The angles between these column vectors are $\arccos \eta$, $\arccos \theta$ and $\arccos \kappa$. Now the formula follows from the cosine theorem of spherical trigonometry (see \cite{Bronstein}, no.3.191), Girard{\textquotesingle}s formula for the area of a spherical triangle (see \cite{Bronstein}, no.3.189), the fact that the volume of the corresponding spherical sector is one third of this area and formula \eqref{volume}.

In case of an elliptic cylinder the inequalities $|\eta |,|\theta |<1$ and $|\kappa |\leqslant 1$ must hold according to Theorem \ref{thm_area}. We transform $P$ with
\begin{equation*}
T:=\left(\begin{matrix} 1 & \eta & \kappa \\ 0 & \sqrt{1-\eta ^{2}} & (\theta -\eta \kappa)/\sqrt{1-\eta ^{2}} \\ 0 & 0 & 1 / \sqrt{1-\eta ^{2}}\end{matrix}\right)
\end{equation*}
in order to get $\Psi(P) = \lbrace (x,y,z) \in \mathbb{R}^3 \mid x^{2} + y^{2} = 1 \rbrace$ because of $|S|=0$. For the columns $t_{1},t_{2},t_{3}$ of $T$ the region of $(\Psi(P),o;t_1,t_2,t_3)$ lies between the $x y$-plane and two planes each of which is defined by a linear function of $x,y$ over a certain circle sector. Hence its volume can be computed by double integration of the two functions. Normal vectors of those planes are given by
\begin{equation*}
t_{1}\times t_{3}=\left(0,-\frac{1}{\sqrt{1-\eta ^{2}}},\frac{\theta -\eta \kappa }{\sqrt{1-\eta ^{2}}}\right) , t_{2}\times t_{3}=\left(1,-\frac{\eta }{\sqrt{1-\eta ^{2}}},\frac{\eta \theta -\kappa }{\sqrt{1-\eta ^{2}}}\right).
\end{equation*}
Hence the demanded functions are
\begin{equation*}
(x,y) \mapsto \frac{y}{\theta -\eta \kappa } \text{ and } (x,y) \mapsto \frac{x \sqrt{1-\eta ^{2}} - \eta y}{\kappa -\eta \theta }
\end{equation*}
over the unit circle sector from (angle) 0 to $\arccos \kappa $ and from $\arccos \kappa $ to $\arccos \eta $, respectively. Hereby the first function must be omitted in case of $\kappa =1$, since this implies $\eta =\theta $ under the condition $|S|=0$, hence the first plane is orthogonal to the $x y$-plane. Otherwise its denominator is positive in case of a bounded sector region.

To prove this we assume $\theta \leqslant \eta \kappa $, i.e. the contrary. Because of $|S|=0$ it follows $\theta =\eta \kappa -\sqrt{(1-\eta ^{2})(1-\kappa ^{2})}$. Boundedness in case of an elliptic cylinder means that the orthogonal projections of $t_{1},t_{2},t_{3}$ onto the $x y$ -plane lie in the same half plane. Since the $y$ -coordinates of $t_{2}$ and $t_{3}$ have opposite sign we must have $\eta +\kappa >0$. Because of $\eta \leqslant \kappa <1$ and $|S|=0$ this implies $-1<\eta \leqslant \theta <0<\kappa $ and $\kappa ^{2}>\eta ^{2}$, hence $\theta <\eta \kappa -(1-\eta ^{2})=\eta (\eta +\kappa )-1<-1$, a contradiction to the boundedness according to Theorem \ref{thm_area}.

The denominator of the second function is not zero\footnote{It is even positive, as the following volume formula shows.} because of $|\eta |,|\theta |<1$ and $\left(1-\eta ^{2}\right)\left(1-\theta ^{2}\right)-(\kappa -\eta \theta )^{2}=|S|=0$. Using polar coordinates we obtain the identity
\begin{equation*}
\mu _{3}\left(P^{+}\right)=\frac{1}{3}\left(\frac{1-\kappa }{\theta -\eta \kappa }+\frac{1-\theta }{\kappa -\eta \theta }\right)=\frac{1-\eta }{3}\left(\frac{\kappa +\theta -\eta -1}{(\theta -\eta \kappa )(\kappa -\eta \theta )}\right)
\end{equation*}
in case $\kappa <1$, i.e. $\theta >\eta \kappa $. In case of $\kappa =1$ the identity reads
\begin{equation*}
\mu _{3}\left(P^{+}\right)=\frac{1}{3}\left(\frac{1-\theta }{1-\eta \theta }\right)=\frac{1}{3(1+\eta )}\;.
\end{equation*}
Now the claimed formulas follow from formula \eqref{volume}.

The case of a plane is trivial.

In case of a hyperbolic cylinder we transform with
\begin{equation*}
T:=\left(\begin{matrix}1 & \eta & \kappa \\ 0 & (\eta \kappa -\theta) / \sqrt{\kappa ^{2}-1} & \sqrt{\kappa^{2}-1} \\ 0 & 1 / \sqrt{\kappa^{2}-1} & 0 \end{matrix}\right)
\end{equation*}
in order to get $\Psi(P) = \lbrace (x,y,z) \in \mathbb{R}^3 \mid x^{2} - y^{2} = 1 \rbrace$ because of $|S|=0$. In analogy with the case of an elliptic cylinder we doubly integrate the functions
\begin{equation*}
(x,y) \mapsto \frac{y}{\eta \kappa -\theta } \text{ and } (x,y) \mapsto \frac{x \sqrt{\kappa ^{2}-1} - \kappa y}{\kappa \theta -\eta} ,
\end{equation*}
using the coordinates $(x,y)=(r\cosh a,r\sinh a)$ with
\begin{equation*}
0\leqslant r\leqslant 1 , 0\leqslant a\leqslant \operatorname{arcosh}\eta \textnormal{ and } \operatorname{arcosh}\eta \leqslant a\leqslant \operatorname{arcosh}\kappa ,
\end{equation*}
respectively. Hereby the first function must be omitted in case $\eta =1$. Then we get the desired result in analogy with the case of an elliptic cylinder. The investigation of boundedness is much simpler: In case $\eta ,\theta ,\kappa \geqslant 1$ the sector region is always bounded.

In case $\theta ,\kappa \geqslant 0$ the function
\begin{equation*}
z(x,y):=\sqrt{l^{2}(x,y)+1-q(x,y)}-l(x,y)
\end{equation*}
with linear form $l(x,y):=\kappa x+\theta y$ and quadratic form $q(x,y):=x^{2}+2\eta xy+y^{2}$ is well defined over the planar sector region
\begin{equation*}
Q^{+} = \lbrace (x,y) \in \mathbb{R}^2 \mid q(x,y)\leqslant 1 ; x,y \geqslant 0 \rbrace .
\end{equation*}
It has the properties $z(x,y)\geqslant 0$ and $(x,y,z(x,y))\in P$ for all $(x,y)\in Q^{+}$. Hence it defines the bounding quadric surface of $P^{+}$. Under the additional condition $\eta \geqslant 0$ it follows
\begin{equation*}
\mu _{3}\left(P^{+}\right)=\int _{0}^{1}I(y)\mathit{dy} \text{ for } I(y):=\int_{0}^{t(y)}z(x,y)\mathit{dx}
\end{equation*}
with $t(y):=\sqrt{1+\eta ^{2}y^{2}-y^{2}}-\eta y$. With help of the pendant in $\mathbb{R}^{3}$ (see \cite{AMR}, Ex. 7.3E) of Leibniz{\textquotesingle} sector rule it turns out that
\begin{equation*}
I(y)=\frac{1}{3} \int_{0}^{t(y)}{\frac{\mathit{dx}}{\sqrt{r(x,y)}}} \text{ with } r(x,y):=l^{2}(x,y)+1-q(x,y)\;.
\end{equation*}
Regarding $r$ as a function of $x\in [0,t(y)]$ for fixed $y\in [0,1]$ it is of the quadratic form $r(x) = a x^2 + b x + c$ with
\begin{equation*}
a := \kappa^{2}-1 , b := 2 (\kappa \theta - \eta) y , c:= 1 + (\theta^2 - 1) y^2 .
\end{equation*}
In case $1\leqslant \eta \wedge 1 < \kappa$ the number $4 a c - b^2 = 4\left(\kappa^{2} - 1 + |S|y^2\right) =: D_{\kappa}(y)$ is always positive thus having the same sign like $\kappa^{2} - 1$. This covers all bounded sector regions of a hyperboloid of two sheets. By looking up a standard integral table \cite{Bronstein}, ch.21, no.241 we get the first formula for the hyperboloid. In case $|S|<0 \wedge 0 \leqslant \eta < 1$ we interchange the roles of $\eta$ and $\kappa$ in the integral, such that $\eta^{2}-1$ and $D_{\eta}(y)$ both have negative sign. By looking up the integral table \cite{Bronstein}, ch.21, no.241 again we find the second formula for the hyperboloid of one sheet. The analysis of boundedness in case of a hyperboloid of one sheet is done via transformation to normal form in analogy with the elliptic cylinder. Hereby it is remarkable that in case $\eta < 0$ the sufficient condition $\eta +\kappa > 0$ of boundedness implies the sufficient condition $\theta > \eta \kappa$ of boundedness. Therefore the former condition can be skipped. For $\beta := \chi(Q,o;p_j,p_k) > -2$ and $x := y := 1/\sqrt{2 + \beta}$ holds $x^2 + \beta x y + y^2 = 1$. Hence $p_{j k} = x p_j + y p_k$ is indeed a point of the bounding arc between $p_j$ and $p_k$. We have $\chi(Q,o;p_j,p_{j k}) = \chi(Q,o;p_{j k},p_k) = \sqrt{2 + \beta} > 0$ and $\chi(Q,o;p_{1 2},p_{2 3}) = (1 + \eta + \theta + \kappa)/\sqrt{1 + \eta + \theta + \eta \theta}$ and an analogous formula with $p_{1 3}$ instead of $p_{1 2}$ or $p_{2 3}$. So the assertion about the four sub sectors follows from the fact $1 + \eta + \theta + \kappa > 1 + \eta + (\eta + 1) \kappa > 0$ for $-1 < \eta$ and $\eta \kappa < \theta \le \kappa$.
\end{proof}

\section{Note on angles}
The measure of the angle that comprises two adjacent angles is the sum of the measures of those two angles. This well known fact about circle sectors generalises naturally to arbitrary planar quadric sectors at centre. To be more precise we generalise the concept {\textquotesingle}angle{\textquotesingle}.

\begin{defn}\label{def_angle}
For linearly independent $a,b \in \mathbb{R}^{n}$ and a $c \in \langle a , b \rangle$ linearly independent from $a$ and from $b$ let $Q$ denote the quadric with external centre $o$ that is determined by $a,b,c \in Q$ (see Theorem \ref{thm_sectorCoeff}b)). In case $\delta :=\delta(a,b;c) := \chi(Q,o;a,b) / 2 > -1$ of a bounded sector region we call
\begin{equation*}
\angle(a,b;c):=\left\{\begin{matrix}\arccos \delta \; \text{if} \; \delta \leqslant 1 \\ \operatorname{arcosh}\delta \; \text{if} \; \delta>1 \end{matrix}\right.
\end{equation*}
the \textit{angle} or \textit{angular measure between} $a$ \textit{and} $b$ \textit{with respect to} $c$.
\end{defn}

The coincidence
\begin{equation*}
\angle (a,b;c)= \arccos \frac{a \circ b}{|a| |b|}
\end{equation*}
with the usual angle (defined via inner multiplication $\circ$) in the special case $|a|=|b|=|c|$ and the symmetry relations
\begin{equation*}
\angle(b,a;c)=\angle(a,b;c)=\angle(a,b;-c)
\end{equation*}
are trivial. We obtain some more interesting properties from Theorem \ref{thm_area}.

\begin{cor}\label{cor_angle}
Let $T$ denote a real, symmetric $2 \times 2$-matrix and
\begin{equation*}
Q := \lbrace (x,y) \in \mathbb{R}^2 \mid (x,y) T (x,y)^t = 1 \rbrace .
\end{equation*}
Let $a,b,c \in Q$ pairwise linearly independent.

a) For the area $\Delta$ of the triangle with vertices $o,a,b$ it holds
\begin{equation*}
\delta^2(a,b;c) = 1 - 4 \Delta^2 \det T .
\end{equation*}

b) In the elliptic or bounded hyperbolic case the region's area of $(Q,o;a,b)$ determined by $c$ equals $\angle(a,b;c)/(2 \sqrt{|\det T|})$.

c) For $(x,y) T (x,y)^t = \alpha x^{2}+ \beta xy + \gamma y^{2}$ and $a = (\kappa,\lambda) , b = (\mu,\nu)$ the region's area of sector $(Q,o;a,b)$ is bounded if and only if
\begin{equation*}
\delta := \alpha \kappa \mu + \beta(\kappa \nu + \lambda \mu)/2 +  \gamma \lambda \nu > -1 .
\end{equation*}
Then it equals
\begin{equation*}
\frac{\arccos(\delta)}{\sqrt{4 \alpha \gamma - \beta^2}} \text{ or } \frac{\operatorname{arcosh}(\delta)}{\sqrt{\beta^{2}-4 \alpha \gamma}}
\end{equation*}
in the elliptic or bounded hyperbolic case, respectively.

d) In case $c$ lies between $a$ and $b$, i.e. $c = x a + y b$ with $x,y>0$, it holds
\begin{equation*}
\angle(a,c;b)+\angle(c,b;a)=\angle(a,b;c) \textit{ (sum of angles in a half plane)} .
\end{equation*}

e) In (the elliptic) case of $-c$ lying between $a$ and $b$ it holds
\begin{equation*}
\angle(a,b;c)+\angle(b,c;a)+\angle(c,a;b) = 2\pi \textit{ (trisection of full ellipse)} .
\end{equation*}
\end{cor}
\begin{proof}
a) For the sector matrix $S$ of $(Q,o;a,b)$ with frame affinity $\Phi$ it holds $1 - \delta^2(a,b;c) = \det S = (\det \Phi)^2 \det T = 4 \Delta^2 \det T$ by formula \eqref{transfFormula}.

b) It is well known that $\det T$ vanishes only in the straight line case. Hence the assertion follows from a) and the area formulae of Theorem \ref{thm_area}.

c) The assertion about boundedness follows also from Theorem \ref{thm_area}. The formula for $\delta$ follows from
\begin{equation*}
\left(\begin{matrix} 1 & \delta \\ \delta & 1 \end{matrix}\right) = \left(\begin{matrix} \kappa & \lambda \\ \mu & \nu \end{matrix}\right) \left(\begin{matrix} \alpha & \beta/2 \\ \beta/2 & \gamma \end{matrix}\right) \left(\begin{matrix} \kappa & \mu \\ \lambda & \nu \end{matrix}\right) .
\end{equation*}

d) In the straight line case $\delta =1$ the identity is trivially fulfilled since all involved angles are $\arccos 1=0$. In the other cases the claimed identity follows from b) and the fact that the biggest area is the sum of the other two.

e) In case $\delta \geqslant 1$ the vector $-c$ can not lie between $a$ and $b$ . So it suffices to consider $\delta <1$. It is well known that an ellipse defined by the equation $\alpha x^{2} + \beta x y + \gamma y^{2} = 1$ has area $\pi /\sqrt{4 \alpha \gamma - \beta^2}$. So the assertion follows from b) and the fact that the three sector regions cover the whole ellipse.
\end{proof}

\begin{rem}
With notation as in Remark \ref{rem_determinants} the vectors
\begin{equation*}
a := (a_1 , a_2) , b := (b_1 , b_2) , c := (c_1 , c_2)
\end{equation*}
are pairwise linearly independent because of $\lambda , \mu , \nu \ne 0$. In the elliptic case the questionable angles are $\varphi := \arccos(\delta) , \psi := \arccos(\varepsilon) , \omega := \arccos(\zeta)$. Assume $\varphi + \psi = \pi$. Since $a$ or $-a$ lies between $b$ and $c$ it follows $\omega = \pi$ by Corollary \ref{cor_angle}d)\&e). But that would mean $\zeta = -1$, a contradiction to Remark \ref{rem_determinants}. So we have $\varphi + \psi \ne \pi$ and analogously $\psi + \omega \ne \pi$ and $\varphi + \omega \ne \pi$. If we apply the Corollary to the two remaining cases $\varphi + \psi > \pi$ and $\varphi + \psi < \pi$ we will obtain: Under the basic requirement $|\delta| < 1$ the inequalities
\begin{equation*}
\delta + \varepsilon > 0 , \varepsilon + \zeta > 0, \delta  + \zeta > 0
\end{equation*}
are equivalent. If $|\delta| < 1$ and $\delta + \varepsilon < 0$ we have
\begin{equation*}
\arccos(\delta) + \arccos(\varepsilon) + \arccos(\zeta) = 2 \pi .
\end{equation*}
If $|\delta| < 1$ and $\delta + \varepsilon > 0$ we have
\begin{equation*}
\arccos(\alpha) + \arccos(\beta) = \arccos(\gamma)
\end{equation*}
for every permutation $(\alpha , \beta , \gamma)$ of $(\delta , \varepsilon , \zeta)$ s.t. $\gamma = \min \lbrace \delta , \varepsilon , \zeta \rbrace$. So we can tell the mutual position of the three ellipse points $a,b,c$ ('full-ellipse-trisection' or not) by the sign of the sum of two arbitrary (of the three) sector coefficients.
\end{rem}


\par\bigskip
\centerline{Eingegangen am  15.\,09.\,2017}
\par\bigskip
Hochschule M\"unchen f\"ur\\
Angewandte Wissenschaften\\
Lothstr. 34\\
D-80335 M\"unchen\\
Germany\\
email\,\,{kahl@hm.edu}

\end{document}